\documentclass[12pt]{article}
\def\version{13.2.2019}\def\WHO{Tom}
\def\users{}  %
\def\users{final-layout}  
\usepackage[utf8]{inputenc}
\usepackage[OT1]{fontenc}
\usepackage{amsmath,amsfonts,amssymb,color}
\usepackage{mathrsfs} 
\usepackage{epsfig}
\usepackage{psfrag}
\usepackage{upgreek} 
\usepackage[format=hang, indention=0cm, labelformat=simple, labelsep=space, textformat=simple, justification=justified, singlelinecheck=true, labelfont=bf, textfont=sl,width=1.0\textwidth,figurename=Figure,skip=15pt]{caption}

\usepackage{cite}

\newtheorem{theorem}{Theorem}[section]
\newtheorem{lemma}[theorem]{Lemma}
\newtheorem{definition}[theorem]{Definition}

\newtheorem{proposition}[theorem]{Proposition}

\newtheorem{remark}[theorem]{Remark}

\usepackage{ifthen}
\ifthenelse{\equal{\users}{final-layout}}{}{
\usepackage{fancyhdr}
\pagestyle{fancy}
\headheight=28pt\headwidth=16.5cm
\definecolor{gray}{gray}{0.5}
\rhead{\color{gray}
T.Roub\'\i\v cek\\Implicit discretisation of damage}
\chead{}
\lhead{File: \jobname.tex on which {\color{red}\WHO}'s 
currently working.\\
Version \version, compiled:
\number\day.\number\month.\number\year\ at
\the\hour:\ifnum\minute<10 0\fi\the\minute\ h }
}

\definecolor{labelkey}{rgb}{1.,.2,0.}

\newcount\hour \newcount\minute
\hour=\time
\divide \hour by 60
\minute=\time
\loop \ifnum \minute > 59 \advance \minute by -60 \repeat

\usepackage[normalem]{ulem}
\definecolor{brown}{rgb}{0.5,0,0}
\ifthenelse{\equal{\users}{final-layout}}{

    \newcommand{\DELETE}[1]{}
    
    \newcommand{\COMMENT}[1]{}
    
    \newcommand{\TINY}[1]{}
    \newcommand{\MARGINOTE}[1]{}
}
{

 \newcommand{\DELETE}[1]{{\color{brown}\sout{#1}\color{black}}}
 
 \newcommand{\COMMENT}[1]{{\color{red}\uuline{#1}\color{black}}}
 
 \newcommand{\TINY}[1]{{\tiny#1}}
 \newcommand{\MARGINOTE}[1]{\marginpar{\color{red}\tiny\texttt{#1}}}
}

\renewcommand\dot[1]{\mathchoice
                 {{\buildrel{\hspace*{.1em}\text{\LARGE.}}\over{#1}}}
                 {{\buildrel{\hspace*{.1em}\text{\Large.}}\over{#1}}}
                 {{\buildrel{\hspace*{.1em}\text{\large.}}\over{#1}}}
                 {{\buildrel{\hspace*{.1em}\text{\large.}}\over{#1}}}}
\newcommand\DT{\dot}
\newcommand\DDT[1]{\mathchoice
   {{\buildrel{\hspace*{.1em}\text{\LARGE.\hspace*{-.1em}.}}\over{#1}}}
   {{\buildrel{\hspace*{.1em}\text{\Large.\hspace*{-.1em}.}}\over{#1}}}
   {{\buildrel{\hspace*{.1em}\text{\large.\hspace*{-.1em}.}}\over{#1}}}
   {{\buildrel{\hspace*{.1em}\text{\large.\hspace*{-.1em}.}}\over{#1}}}}
\newcommand{\wt}[1]{\mathchoice
     {\text{\small$\widetilde{\text{\normalsize$#1$}}\hspace*{.03em}$}}
                    {\text{\small$\widetilde{\text{\normalsize$#1$}}$}}
                    {\widetilde{#1\hspace*{-.02em}}\hspace*{.03em}}
                    {\widetilde{#1}}}

\newcommand{\lineunder}[2]{\LU{\begin{array}[t]{c}\underbrace{#1}\vspace*{.5em}\end{array}}{\mbox{\footnotesize\rm #2}}}

\newcommand{\LU}[2]{\begin{array}[t]{c}#1\vspace*{-1em}\\_{#2}\end{array}}

\def\R{{\mathbb R}}
\newcommand\eps{\varepsilon}
\newcommand\bbC{\mathbb C}
\newcommand\bbD{\mathbb D}

\newcommand\bbI{\mathbb I}

\renewcommand\d{\mathrm d}

\newcommand{\NEW}[1]{{\color{black}#1}}

\begin{document}
\begin{sloppypar}
\def\taur{\chi_{_{\rm R}}}

\begin{center}
\LARGE\bfseries
Coupled time discretisation of dynamic damage models at small strains\footnote{{\tt This is an author accepted version published on line on April 2019
    in {\it IMA Journal of Numerical Analysis},
    DOI 10.1093/imanum/drz014,
     later as  (2020) {\bf 40}, 1772--1791, and
     placed on arXiv in accord with the ``Standard Licence'' of the Oxford University Press.}}$^{,}\!\!$
    \footnote{{\tt  URL https://academic.oup.com/imajna/article-abstract/40/3/1772/5431182, }}

\end{center}

\bigskip

\begin{center}
\bfseries
Tom\'{a}\v{s} Roub\'\i\v{c}ek\footnote{
Institute of Thermomechanics of the Czech Acad. Sci., Dolej\v skova~5,
CZ--182~08 Praha 8, Czech Republic, email: {\tt tomas.roubicek@mff.cuni.cz}}
\end{center}

\begin{abstract}
The dynamic damage model in viscoelastic materials in Kelvin-Voigt rheology 
is discretised by a scheme which is coupled, suppresses spurious numerical
attenuation during vibrations, and has a variational structure with a convex
potential for small time-steps. In addition, this discretisation is
numerically stable and convergent for the time step going to zero. When
combined with the FEM spatial discretisation, it leads to an implementable
scheme and to that iterative solvers (e.g.\ the Newton-Raphson)
used for the nonlinear algebraic systems at each
time level have guaranteed global convergence. Models which are
computationally used in some engineering simulations
in a non-reliable way are thus stabilized and theoretically justified in this
viscoelastic rheology. In particular, this model and algorithm can be used
in a reliable way for a dynamic fracture in the usual phase-field
approximation.

\medskip


\medskip

\end{abstract}

\section{Introduction}\label{sec-intro}

Damage (possibly interpreted as a phase-field approximation of fracture)
is an important phenomenon in continuum mechanics of solids and many models
have been devised in engineering and partly also analyzed in mathematical
literature, cf.\ e.g.\ \cite{KruRou18MMCM,Roub??MDDP} and references therein.
Mostly, a scalar-valued damage field is considered. 
Historically, this concept dates back to Kachanov \cite{Kach58TRPD}, and has
been widely used both at finite and at small strains, cf.\ e.g.\ the monographs
\cite{Frem02NST,Kach90ICDM,Maug92TPF,ShSoTe04MAQC,SoHaSh06AACP}. 

In agreement with experiments, one other aspect is often built into damage 
models, namely an internal length scale through the \emph{gradient of damage}, 
cf.~\cite{BazJir02NIFP,Frem02NST},
and, in some situations, also \emph{inertial effects}. 
The former aspect expresses certain
nonlocality in the sense that damage of a particular spot is to some extent 
influenced by its surrounding, leading to possible hardening or 
softening-like effects and, by introducing a certain internal length scale,
eventually prevents damage microstructure development. 
The latter aspect is often important, in particular because damage
may evolve very fast and elastic waves or vibrations may vitally accompany
the damage process itself, and also the finite speed of propagation of
information naturally avoids nonphysical long-range interactions
sometimes involved in quasistatic damage models. This happens
in geophysics during ruptures of lithospheric faults when seismic waves are
emitted or in material testing in engineering. 
In some other situations, the external load may be extremely fast (e.g.\ during
explosions or impacts of projectiles  or asteroids) so that
dynamical effects may play a vital role.

Although the split (called also fractional-step, or staggered) schemes are
robust and even may conserve energy, cf.\ also
\cite{MiWeHo10PFMR,MolGra17AIRS,Roub17ECTD,RouPan17ECTD}
or \cite[Sect.\,5.1-5.2]{MieRou15RIST},
the fully-coupled scheme are often used in engineering 
inspite of no guaranty for numerical stability and convergence,
or even for a global convergence of iterative solvers
which must be necessarily used for
(generally) nonconvex 
incremental problems arising at each time level. To the last point,
typically the Newton-Raphson or the alternate minimization algorithm (AMA)
are used. Especially the former option is explicitly reported in literature
as unreliable in the context of damage or phase-field crack models,
cf.\ \cite{BVSH12PFDD,GADS12NAMQ,NeAlDo09NMCB,RoyDod01SDCG,SoPaBu06SCPA,SoPaBu06BCZM},
although one can even find an opposite (not correct in general) belief and a
(rather particular) experience that a fully coupled monolithic solution can
even prevent numerical instabilities or converges, cf.\ e.g.\
\cite{PaCoRe18FSGA,PagRei17RPCI}.
As a matter of fact, the Newton-Raphson algorithm, if converges, can find only
a general critical points which does not need to be physically
relevant and cannot yield apriori estimates and thus numerical
stability of the discretisation scheme.

\NEW{Let us only mention that} the AMA procedure works
\NEW{differently}
because it monotonically increases energy and, if converges, it
founds critical points which minimize the energy at least separately at
the displacement and at the damage variables,
cf.\ \cite{Bour07NIVF} for a theoretical analysis. The relation
of AMA and the Newton-Raphson algorithm is discussed in \cite{BuOrSu10AFEA}.
Sometimes, the variant of AMA which updates the damage constraints within
each iteration is used, which is then basically the staggered scheme, 
cf.\ \cite{KneNeg17CMS}. Cf.\ also \cite{FarMau17LNSV} for a comparison
and some other variants.

Sometimes, attempt for stabilization of such unreliable schemes has
intuitively been done by inventing some viscosity into the damage flow rule
\cite{KnRaBr01CPFA,MCMD07CDMC,MiWeHo10TCPF,WilKno09SDFO},
although even this is not reliable, cf.\ Remark~\ref{rem-visco-dam}.
For some others, even more phenomenological and theoretically
unjustified attempts to cope with failing iterative solvers, see e.g.\ 
\cite{NeAlDo09NMCB}. 

The goal of this article is to investigate a true physically-motivated
stabilization by means of
the viscosity in the bulk, i.e.\ by inventing Kelvin-Voigt rheology
instead of a purely elastic model, or a less robust stabilization
by mere inertia in a purely elastic model, cf.~Remark~\ref{rem-purely-elastic}.
This viscosity has a physical motivation in particular during fast
crack growth or fast loading, cf.\ the discussion in \cite[\$8.2]{Freu98DFM}.
The fully-coupled schemes thus become amenable for a-priori estimates 
(numerical stability) and their convergence can be shown, too.
In addition, 
we keep the variational structure of the incremental problems and,
for small time step, the guaranty of global convergence of iterative
numerical algorithms.

It should be emphasized that {\it fully coupled} (also, in engineering
literature, sometimes called 
{\it monolithic}) schemes do not rely on component-wise convex (or
even quadratic) structure
and, if working, they can be used for more general problems
without any algorithmic changes, in contrast to split schemes that
typically take a benefit from implementation of quadratic-programming.
In particular, it concerns the very natural phenomenon (not considered
here) that damage may act very differently on compression than on tension
so that the stored energy is not quadratic in terms of strains.
Even stored energies nonconvex in terms of strains 
are considered in some applications in rupturing-rock mechanics, cf.\ 
\cite{HaLyBZ11ESED,LyHaBZ11NLVE} and references therein.
This is a certain advantage of the monolithic schemes, beside their
more straightforward implementation in some software packages as e.g.\ Abaqus,
cf.\ e.g.\ \cite{LLMZ16AIMS,MSJA15AIPF}.

The plan of this paper is following:
The model is formulated in Sect.~\ref{sec-model} by specifying the governing 
energies and by using the Hamilton variational principle. 
Its time discretisation is devised in Sect.~\ref{sec-disc} and its basic 
properties are proved, namely the variational structure, existence of 
discrete solutions, and convexity of the governing potentials for sufficiently 
small time steps. In Sect.~\ref{sec-anal}, the numerical stability
analysis (i.e.\ a-priori estimates) of the time-discrete model
is performed, and its 
convergence towards the original continuous model is proved.
Eventually, the paper ends with several remarks in Sect.~\ref{sec-rem}
outlining the phase-field fracture, and various modifications and further
usage.

\section{The model of damage at small strains}\label{sec-model}

We deal with a relatively simple model of damage in Kelvin-Voigt 
viscoelastic solid without influencing the viscous part. 
We consider also inertial effects, so that we can capture also vibrations or
waves that can be generated during fast damage or rupture. 
Yet, our results do not rely on the presence of inertia so that the 
quasistatic variant is covered, too.

The independent variables of the model are the displacement 
$u:\Omega\to\R^d$ and the damage $\alpha:\Omega\to[0,1]$, 
both varying in time, with $\Omega\subset\R^d$ being a fixed 
bounded domain in $\R^d$ with a Lipschitz boundary $\Gamma$, $d=2,3$.
%
The viscoelastic material is considered with 
the viscous response undergoing (quite naturally) damage in
the same way as the elastic response.
The mass density $\varrho$ which is responsible for inertia is naturally
independent of damage because the mass is not destroyed by 
damaging of interatomic links. 

\def\zetai{\zeta}
\def\zetav{\zeta_{\text{\tiny\rm VI}}^{}}
\def\xii{\xi_{\text{\tiny\rm IN}}^{}}
\def\xiv{\xi_{\text{\tiny\rm VI}}^{}}
\def\zetaistar{\zeta_{\text{\tiny\rm IN}}^*}
\def\zetaiprime{\zeta'}
\def\zetavprime{\zeta_{\text{\tiny\rm VI}}'}


The rational-mechanical approach builds models on the base of energies.
Here, we want to use a model governed by the energies: 
\begin{subequations}\label{energies}\begin{align}\label{energies-E}
&\text{stored energy:}&&\mathscr{E}(u,\alpha):=
\int_\Omega\frac12\bbC(\alpha)e(u){:}e(u)-\phi(\alpha)
+\frac\kappa p|\nabla\alpha|^p\,\d x,
\\\label{energies-D}
&\text{dissipation potential:}&&\mathscr{D}(\alpha;v,\DT\alpha):=
\int_\Omega\frac12\bbD(\alpha)e(v){:}e(v)+\zetai(\DT\alpha)\,\d x,
\\\label{energies-T}
&\text{kinetic energy:}&&\mathscr{T}(v):=\int_\Omega\frac\varrho2|v|^2\,\d x,
\\\label{energies-F}
&\text{external loading:}&&\mathscr{F}(t,u):=
\int_\Omega f(t)\!\cdot\!u\,\d x+\int_\Gamma g(t)\!\cdot\!u\,\d S,
\end{align}\end{subequations}
where $\varrho>0$ is the mentioned mass density, the variable $v$ is the
placeholder for velocity $\DT u$, i.e.\ the time derivative of $u$, 
$e(u)=\frac12\nabla u^\top\!+\frac12\nabla u$ the small-strain tensor,
$\phi=\phi(\alpha)$ is the stored
energy of damage reflecting the phenomenon that, on a microscopical level,
damage means some microcracks or microvoids which increase micro-surfaces 
inside material and thus its internal stored energy. Thus, physically,
$\phi$ is increasing or at least nondecreasing. The negative sign in
\eqref{energies-E} is thus related with the convention that
$\alpha=1$ means no damage while $\alpha=0$ means maximal damage.

Based on the energies \eqref{energies}, 
the evolution is governed by the Hamilton variational principle generalized 
for nonconservative systems. More specifically,
defining the Lagrangian $\mathscr{L}(v,u,\alpha)$ and the nonconservative
dissipative force $\mathscr{F}_{\rm diss}^{}(v,\alpha,\DT\alpha):=
\mathscr{D}_{(v,\DT\alpha)}'(\alpha;v,\DT\alpha)$ by
\begin{subequations}\begin{align}
&\mathscr{L}(t,v,u,\alpha)=
\mathscr{T}(v)-\mathscr{E}(u,\alpha)-\mathscr{F}(t,u),
\\&\langle\mathscr{F}_{\rm diss}^{}(v,\tilde\alpha,\DT{\tilde\alpha}),(u,\alpha)\rangle=
\int_\Omega\bbD(\tilde\alpha)e(v){:}e(u)
+\zetaiprime(\DT{\tilde\alpha})\alpha\,\d x.
\end{align}\end{subequations}
The evolution on a fixed time interval $I=[0,T]$ will be a trajectory 
$t\mapsto(u(t),\alpha(t))$ with $t\in I$ forming a critical point of
the functional 
\begin{align}\label{Ham}
(u,\alpha)\mapsto\int_0^T\mathscr{L}(t,\DT u,u,\alpha)
-
\langle\mathscr{F}_{\rm diss}^{}(v,\tilde\alpha,\DT{\tilde\alpha}),(u,\alpha)\rangle\,\d t
\end{align}
provided $v=\DT u$ and $\tilde\alpha=\alpha$.
For $\mathscr{F}_{\rm diss}^{}=0$, this is known as the Hamilton variational 
principle for conservative systems. The extension \eqref{Ham} for nonconservative
dissipative force $\mathscr{F}_{\rm diss}^{}\ne0$ is a bit formal, cf.\ 
\cite{Bedf85HPCM}.

From \eqref{Ham}, we obtain the system of partial differential
equations (or inclusions): 
\begin{subequations}\label{KV-damage-}\begin{align}\label{KV-damage-u-}
&\varrho\DDT{u}
-{\rm div}\big(
\bbC(\alpha)e(u){+}\bbD(\alpha)e(\DT u)\big)=f
&&\text{in }\ Q,\\\label{KV-damage-z-}
&\partial\zetai(\DT\alpha)+
\frac12\bbC'(\alpha)e(u){:}e(u)
-{\rm div}(\kappa|\nabla\alpha|^{p-2}\nabla\alpha)
\ni\phi'(\alpha)&&\text{in }\ Q,
\end{align}\end{subequations}
where $Q:=I\times\Omega$, accompanied by the boundary conditions 
$(\bbC(\alpha)e(u){+}\bbD(\alpha)e(\DT u))\vec{n}=g$ 
and $\nabla\alpha{\cdot}\vec{n}=0$ on $\Sigma:=I\times\Gamma$,
with $\vec{n}$ denoting the outward unit normal to $\Gamma$.
The notation ``$\partial\zetai$'' stands for the usual convex 
subdifferential, i.e.\ formally \eqref{KV-damage-z-} means
the inequality
\begin{align}
\zetai(z)+
\Big(\frac12\bbC'(\alpha)e(u){:}e(u)
-{\rm div}(\kappa|\nabla\alpha|^{p-2}\nabla\alpha)
-\phi'(\alpha)\Big)(z-\DT\alpha)\ge\zetai(\DT\alpha)
\label{VI-ineq}\end{align}
to be valid on $Q$ for all $z\in\R$.

Throughout this article, we will accept the modelling assumptions
\begin{subequations}\label{model-ass}\begin{align}\label{model-ass1}
&\bbC'(0)=0,\ \ \phi'(0)\ge0,\ \ 
\text{ and }\ \ \zetai(\DT\alpha)=+\infty\ \text{ for }\ 
\DT\alpha>0,
\\&\label{model-ass2}
\bbD(\alpha)=\bbD_0+\taur\bbC(\alpha)\ \ \ 
\text{with $\ \taur>0\ $ a fixed relaxation time}.
\end{align}\end{subequations}
The assumption \eqref{model-ass1} ensure that the evolution 
of damage $\alpha$ is {\it unidirectional} in the sense
that healing $\DT\alpha>0$ is not allowed, and that $\alpha$ stays 
valued in $[0,1]$ within the evolution provided it is such at the
initial time, so that these constraints $0\le\alpha\le1$ do not need to be
considered explicitly, which simplifies the formulation as well as the
analysis. Like in \cite{LaOrSu10ESRM}, \eqref{model-ass2} represents a bit 
special but not much application-restricting case that viscous dissipative 
processes are of the same character as the elastic storage.
Some results (in fact, all except Proposition~\ref{prop3})
would hold even without assuming \eqref{model-ass2}.
In fact, \eqref{model-ass2} is only needed for proving 
a strong convergence of the strains without relying 
on the energy conservation in the damage flow rule, which is
not guaranteed in our rate-dependent uni-directional damage model
(in contrast to the rate-independent model as in \cite{LaOrSu10ESRM} or
\cite[Sect.\,5.2.5]{MieRou15RIST}).  

The energetics of the system \eqref{KV-damage-} is revealed when
testing \eqref{KV-damage-u-} by $\DT u$ and \eqref{KV-damage-u-} by
$\DT\alpha$. This gives, at least formally \NEW{(i.e.\ it the solution
  would be smooth enough so that ${\rm div}(|\nabla\alpha|^{p-2}\nabla\alpha)$
  would be in $L^2(Q)$ and thus in duality with $\DT\alpha$)}, the energy
balance
\begin{align}\label{energetics}
\mathscr{T}(\DT u(t))+\mathscr{E}(u(t),\alpha(t))
+\!\int_0^t\!\!\varXi(\alpha;\DT u,\DT\alpha)\,\d t=
\mathscr{T}(v_0)+\mathscr{E}(u_0,\alpha_0)
+\!\int_0^t\!\!\mathscr{F}(t,\DT u)\,\d t,
\end{align}
where the dissipation rate is
\begin{align}\label{energetics+}
  \varXi(\alpha;v,\DT\alpha)=
v\cdot\partial_v^{}\mathscr{D}(\alpha;v)+
\DT\alpha\partial_{\DT\alpha}\mathscr{D}(\DT\alpha)
=\int_\Omega\bbD(\alpha)e(v):e(v)+\DT\alpha\partial\zeta(\DT\alpha)\,\d x\,.
\end{align}
Typically, $\zeta$ is smooth except only at $\DT\alpha=0$, so that
$\DT\alpha\partial\zeta(\DT\alpha)$ is well defined as a single-valued
function on $Q$.

In terms of the variable $v$ for velocity, instead of the 
2nd-order force-equilibrium equation \eqref{KV-damage-u-} involving the 
inertial force $\varrho\DDT{u}$, we rewrite \eqref{KV-damage-} in the form 
of a 1st-order system
\begin{subequations}\label{KV-damage}\begin{align}
&\DT{u}=v&&\text{in }\ Q,\label{KV-damage-v}
\\&\label{KV-damage-u}
\varrho\DT{v}
-{\rm div}\big(\bbD_0v+
\bbC(\alpha)e(u{+}\taur v)
\big)=f
&&\text{in }\ Q,\\\label{KV-damage-z}
&\partial\zetai(\DT\alpha)+
\frac12\bbC'(\alpha)e(u){:}e(u)
-{\rm div}(\kappa|\nabla\alpha|^{p-2}\nabla\alpha)
\ni 
\phi'(\alpha)
&&\text{in }\ Q.
\intertext{We consider the initial-boundary value problem}
&\big(\bbD_0v+
\bbC(\alpha)e(u{+}\taur v)\big)\vec{n}=g
\ \ \text{ and }\ \ \nabla\alpha{\cdot}\vec{n}=0&&\text{on }\ \Sigma,
\\& u(0)=u_0,\ \ \ \ \ \DT u(0)=v_0,\ \ \ \ \ \alpha(0)=\alpha_0&&\text{on }\ \Omega.
\end{align}\end{subequations}

A very similar model has already been analysed by C.J.\,Larsen, C.\,Ortner,
and E.\,S\"{u}li by using the semi-implicit (fractional-step)
time-discretisation in \cite{LaOrSu10ESRM} for the rate-independent
unidirectional damage, cf.\ also
\cite[Sect.5.1.1 and 5.2.5]{MieRou15RIST}. 
%
In fact, for a special choice of the damage energy $\phi$, 
and for $p=2$, our model nearly equals to that one in \cite{LaOrSu10ESRM}
except that we consider a 
rate-dependent damage evolution, in contrast to \cite{LaOrSu10ESRM} where
a rate-independent damage is considered. Yet, to gain the
variational structure of our monolithic time disretisation,
we will make a semi-implicit discretisation of the viscous term and
will need $p>d$ in the proof of convergence in Proposition~\ref{prop3}.
\NEW{Of course, the choice $p=2$ is more straightforward and common
  in engineering computations. Here, let us only remark that a compromise model
  replacing 
  $|\nabla\alpha|^{p-2}$ in \eqref{KV-damage} by $1+\eps|\nabla\alpha|^{p-2}$
  with some (pressumably small) $\eps>0$ and again $p>d$
would work, too.}
  
The rate-independent models
represent a certain asymptotical abstraction for slow external loading 
activating internal processes which may evolve much faster, and can
sometimes be simpler than their rate-dependent variants, cf.\ 
\cite[Example 1.2.7]{MieRou15RIST} for a damage problem.
Even the rate-dependent model allows for analytical existence result,
cf.\ \cite{KruRou18MMCM,Roub??MDDP} for $p=2$; the existence results 
presented here are only new for $p\ne2$. 
Yet, our focus is on numerical approximation scheme and on its stability
and convergence. In these aspects, there are several differences 
comparing to \cite{LaOrSu10ESRM} beside rate-dependent variant of the
damage-flow rule: the fully coupled scheme is used here (while
\cite{LaOrSu10ESRM} uses the split, staggered scheme) and
the inertial term is here discretised in an energy-conservative way to
suppress spurious numerical attenuation during vibrations which
would otherwise effectively kill computational simulations of such
problems.

\section{Coupled implicit time discretisation}\label{sec-disc}

Although for some $\bbC(\cdot)$ and $\bbD(\cdot)$ the functionals
\eqref{energies-E} and \eqref{energies-D} can be convex, the damage models   
inevitably need to work with nonconvex energies. This
nonconvexity is responsible for the desired phenomenon of sudden
damage/rupture.
Some engineering algorithms intentionally want to keep this nonconvexity
even on the discrete level, exploiting also that software packages as e.g.\
Abaqus accent rather coupled, so-called 
{\it monolithic scheme}, cf.\ e.g.\ \cite{LLMZ16AIMS}.

We apply here the coupled implicit time discretisation to the 
system \eqref{KV-damage} rather than \eqref{KV-damage-} in a 
particular way. Namely: 
\begin{itemize}
\item we discretise the inertial part by the mid-point (Crank-Nicolson)
formula rather than the backward Euler one in order to reduce unwanted 
numerical attenuation which would otherwise be exhibited by the
implicit 2nd-order scheme, 
and
\item
we intentionally use the fully implicit backward-Euler 
formula for the stored-energy terms 
to be as close as \NEW{possible to} the usual engineering implementation,
\item
we  use a semi-implicit (but not the fully implicit backward-Euler) 
formula for the viscous stress $\bbD(\alpha)e(v)$ in order to keep the 
variational structure of the incremental problems, cf.~Proposition~\ref{prop} 
below.
\end{itemize}
The result\NEW{ing iterative} 
coupled boundary-value problems here are:
\begin{subequations}\label{IBVP-damage-small-disc}
\begin{align}
&\frac{u_\tau^k{-}u_\tau^{k-1}}\tau=v_\tau^{k-1/2}\qquad\ \text{ with }\ \ 
v_\tau^{k-1/2}:=\frac{v_\tau^k{+}v_\tau^{k-1}}2,\\\nonumber
&\varrho\frac{v_\tau^k{-}v_\tau^{k-1}}\tau-{\rm div}\,
\Big(
\bbC(\alpha_\tau^{k})e(u_\tau^{k})+\bbD(\alpha_\tau^{k-1})e(v_\tau^{k-1/2})
\Big)
=f_\tau^k
\\&\qquad\qquad\qquad\qquad\qquad\text{ with }\ \ 
f_\tau^k:=\int_{(k-1)\tau}^{k\tau}\!\!f(t)\,\d t
,\ \text{ and }
\label{IBVP-damage-disc-1}
\\
&
\partial\zetai
\Big(\frac{\alpha_\tau^k{-}\alpha_\tau^{k-1}}\tau\Big)
+
\frac12\bbC'(\alpha_\tau^{k})e(u_\tau^{k}){:}e(u_\tau^{k})
-{\rm div}(\kappa|\nabla\alpha_\tau^{k}|^{p-2}\nabla\alpha_\tau^{k})\ni 
\phi'(\alpha_\tau^{k})
\label{IBVP-damage-disc-2}
\end{align}
\end{subequations}
considered on $\Omega$ while completed with the corresponding boundary 
conditions
\begin{subequations}\label{IBVP-damage-small-disc-BC}
\begin{align}\label{IBVP-damage-small-disc-BC1}
&\Big(
\bbC(\alpha_\tau^{k})e(u_\tau^{k})
+\bbD(\alpha_\tau^{k-1})e(v_\tau^{k-1/2})
\Big){\cdot}\vec{n}=g_\tau^k\ \ \text{ and}
\\&\kappa\nabla\alpha_\tau^{k}{\cdot}\vec{n}=0,\ \
\text{ where }\ \ 
g_\tau^k:=\int_{(k-1)\tau}^{k\tau}\!\!g(t)\,\d t.
\label{IBVP-damage-small-disc-BC2}
\end{align}\end{subequations}
It is to be solved 
\NEW{iteratively} for $k=1,...,T/\tau$ with 
\begin{align}
u_\tau^0=u_0,\ \ \ \ \ \ 
v_\tau^0=v_0,
\ \ \ \ \ \ \alpha_\tau^0=\alpha_0.
\end{align}

The important property of the underlying energy functionals is convexity.
In damage-type problems, the mentioned nonconvexity of the
stored energy is suppressed in the time-discrete problems 
due to the quadratic viscosity potentials, and it suffices to guarantee  
only some semiconvexity of the stored energy,
cf.\ \cite[Remark 8.24]{Roub13NPDE} 
or for damage-type problems \cite{Roub09RIPV}. We should emphasize that 
the attribute of semiconvexity is not automatic and e.g.\ for $\bbC(\cdot)$
affine, i.e.\ for $(e,\alpha)\mapsto\alpha\bbC'e{:}e$, it does not hold
cf.\ \eqref{det-hessian+} below for $\gamma''=0$ and realize that
even considering also the term $\phi$ cannot help for $|e|\to\infty$.
Anyhow, here fortunately, we are able to state even a bit strengthened version
of the semiconvexity, using the convexification by only the $e$-component but
not $\alpha$. This will enable to devise the robust scheme by using the
viscosity only in strains but not in $\alpha$, although we will use (and
later need) the rate-dependent flow rule for $\alpha$ for simplifying some
analytical arguments. 

Let us first summarize the assumptions needed partly here and especially
later:
\begin{subequations}\label{ass}
  \begin{align}\nonumber
 &\text{$\bbC(\alpha):=
  \gamma(\alpha)\bbC_1$  with
  $\bbC_1$ positive definite and}
    \\&\qquad\qquad\qquad\quad\ \label{ass1}
  \text{with $\gamma:\R\to\R^+$ smooth, 
  positive, and strictly convex,}
 \\\label{ass1+}
 & \text{$\bbD_0(\cdot)$ in \eqref{model-ass2} symmetric positive definite,}
 \\\label{ass1++}
   & \text{$\phi:[0,1]\to\R$ continuously differentiable \NEW{and concave},}
  \\\nonumber
&\text{$\zetai:\R\to\R\cup\{+\infty\}$ convex lower semicontinuous and}
   \\&\qquad\qquad\qquad\qquad\ \ \,\label{ass1+++}
  \text{coercive in the sense that $\inf_{v\ne0}\zetai(v)/|v|^2>0$,}
   \\\label{ass2}
& \text{$\varrho\in L^\infty(\Omega),\ \ $ess\,inf$_\Omega\varrho(\cdot)>0$,\ \ $f\in L^2(Q;\R^d)$, $g\in L^2(\Sigma;\R^d)$,}
\\\label{ass3}
& \text{
$u_0\in H^1(\Omega;\R^d)$, $v_0\in L^2(\Omega;\R^d)$, and
  $\alpha_0\in W^{1,p}(\Omega)$.}
\end{align}\end{subequations}

\begin{lemma}[Strenghtened semiconvexity of the stored energy.]\label{lem}
  Let \eqref{ass1} hold. Then the function 
$(e,\alpha)\mapsto
\frac12\bbC(\alpha)e{:}e+\frac12K|e|^2
$ 
is convex
provided $K$ large enough.
\end{lemma}

\noindent{\it Proof.}
The Hessian of the function in question, i.e.\ 
\begin{align}
  \bigg(\!\!\begin{array}{ccc} 
    \gamma(\alpha)\bbC_1+K\bbI\!\!\!&,&\!\!\!\gamma'(\alpha)\bbC_1e\\
\gamma'(\alpha)\bbC_1e\!\!\!&,&\!\!\!\frac12\gamma''(\alpha)\bbC_1e{:} e
\end{array}\!\!\bigg)&=
  \bigg(\!\!\begin{array}{ccc} 
    \gamma(\alpha)\bbC_1\!\!\!&,&\!\!\!0\\0\!\!\!&,&\!\!\!0
\end{array}\!\!\bigg)
  +\bigg(\!\!\begin{array}{ccc}K\bbI\!\!\!&,&\!\!\!\gamma'(\alpha)\bbC_1e\\
\gamma'(\alpha)\bbC_1e\!\!\!&,&\!\!\!\frac12\gamma''(\alpha)\bbC_1e{:} e
\end{array}\!\!\bigg)
\,,
\label{hessian}\end{align}
is to be positive semidefinite for sufficiently big $K$.
The former matrix on the right-hand of \eqref{hessian}
is  surely positive semidefinite. Therefore we are to prove the positive
semidefiniteness of the latter one. \NEW{For any
  $(\wt e,\wt\alpha)\in\R_\text{sym}^{d\times d}\times\R$, we can estimate}
\begin{align}\nonumber
 \NEW{\bigg(\!\!\begin{array}{c}\wt e\\\wt\alpha\end{array}\!\!\bigg)^\top
    \bigg(\!\!\begin{array}{ccc}K\bbI\!\!\!&,&\!\!\!\gamma'(\alpha)\bbC_1e\\
\gamma'(\alpha)\bbC_1e\!\!\!&,&\!\!\!\frac12\gamma''(\alpha)\bbC_1e{:} e
    \end{array}\!\!\bigg)
    \bigg(\!\!\begin{array}{c}\wt e\\\wt\alpha\end{array}\!\!\bigg)}
    &\NEW{=K|\wt e|^2+\frac12\wt\alpha^2\gamma''(\alpha)\bbC_1e{:} e
    +2\wt\alpha\gamma'(\alpha)\bbC_1e{:}\wt e}
    \\&\nonumber\NEW{\ge\wt\alpha^2\Big(\frac12\gamma''(\alpha)\bbC_1e{:} e
    -\frac1K\gamma'(\alpha)^2|\bbC_1e|^2\Big)}
\\&\NEW{\ge\wt\alpha^2}\Big(\frac{\gamma''(\alpha)}{2|\bbC_1^{-1}|}
-\frac1K\gamma'(\alpha)^2|\bbC_1|^2\Big)|e|^2
  \ge0
\,,
\label{det-hessian+}\end{align}
where $1/|\bbC_1^{-1}|$ is the positive-definiteness constant of $\bbC_1$, i.e.\
$1/|\bbC_1^{-1}|=\min_{|e|=1}\bbC_1e{:}e$. It reveals that, for the last
inequality in \eqref{det-hessian+}, it suffices to take
\begin{align}\label{K-large}
  K\ge
  2
  |\bbC_1|^2|\bbC_1^{-1}|
  \frac{\max\gamma'([0,1])^2}{\min\gamma''([0,1])}\,.
\end{align}
$\hfill\Box$

\medskip

The important attribute behind the scheme 
\eqref{IBVP-damage-small-disc}--\eqref{IBVP-damage-small-disc-BC}
is its variational structure, holding even for any time steps $\tau>0$ 
without the mentioned convex structure,  
cf.\ also \cite{Roub17ECTD,RouPan17ECTD} for the form of \eqref{potential}
if the stored energy were convex:

\begin{proposition}[Variational structure.]\label{prop}
  Let \eqref{model-ass} and {\rm(\ref{ass}a-e)} be valid,
  $u_\tau^{k-1}\in H^1(\Omega;\R^d)$, $v_\tau^{k-1}\in L^2(\Omega;\R^d)$,
  $\alpha_\tau^{k-1}\in L^\infty(\Omega)$, $0\le\alpha_\tau^{k-1}\le1$,
  and let $p>1$. Then  the functional
\begin{align}\nonumber\\[-2.2em]\nonumber
(u,\alpha)\mapsto&\int_\Omega
\frac{\varrho}{2\tau}\Big|\frac{u{-}u_\tau^{k-1}}\tau{-}v_\tau^{k-1}\Big|^2\!
+\frac12\bbC(\alpha)e(u){:}e(u)-\phi(\alpha)+\frac\kappa p|\nabla\alpha|^p
-f_\tau^k{\cdot}u
\\[-.2em]&
\
+\frac{1}{2\tau}\bbD(\alpha_\tau^{k-1})e(u{-}u_\tau^{k-1}){:}e(u{-}u_\tau^{k-1})
+\tau\zetai\Big(\frac{\alpha{-}\alpha_\tau^{k-1}}\tau\Big)
\,\d x-\int_\Gamma g_\tau^k{\cdot}u\,\d S
\label{potential}\end{align}
is weakly lower semicontinuous on $H^1(\Omega;\R^d)\times H^1(\Omega)$.
For any $(u_\tau^k,v_\tau^k,\alpha_\tau^k)\in H^1(\Omega;\R^d)\times 
L^2(\Omega;\R^d)\times W^{1,p}(\Omega)$ solving (in the usual weak sense)
the boundary value problem 
\eqref{IBVP-damage-small-disc}--\eqref{IBVP-damage-small-disc-BC},
the couple  $(u_\tau^k,\alpha_\tau^k)$ is a critical point of this functional. 
Also, conversely, any critical point $(u,\alpha)$ of \eqref{potential}
gives a weak solution  $(u_\tau^k,v_\tau^k,\alpha_\tau^k)$ 
to \eqref{IBVP-damage-small-disc}--\eqref{IBVP-damage-small-disc-BC} 
when putting 
\begin{align}
u_\tau^k=u,\qquad
v_\tau^k=\frac2\tau(u_\tau^k{-}u_\tau^{k-1})-v_\tau^{k-1},\qquad\alpha_\tau^k=\alpha.
\label{potential+}
\end{align}
In particular, such solutions do exist. Moreover, 
the potential \eqref{potential} is strictly convex provided
$\tau$ is sufficiently small, namely   
\begin{align}\nonumber\\[-1.8em]\label{tau-small}
  \ \ \ \ \ \ \ \ \tau\le\tau_0=\frac{\min\gamma''([0,1])}{2|\bbD_0^{-1}
    |\,|\bbC_1|^2|\bbC_1^{-1}|\max\gamma'([0,1])^2}\quad
  \text{with $\bbD_0:=\bbD(0)$ as in \eqref{model-ass2}}\,.
\end{align}
\end{proposition}


\noindent{\it Proof.} 
The system \eqref{IBVP-damage-small-disc} does not satisfy the usual symmetry 
condition and thus does not have any potential, but eliminating $v_\tau^k$ by 
substituting $v_\tau^k$ as in \eqref{potential+}, we obtain a potential for
the couple $(u_\tau^k,\alpha_\tau^k)$, and it is not difficult to 
identify its form as \eqref{potential}.

For the convexity, we use Lemma~\ref{lem} for the $\bbC$-term, which is
the only nonconvex term in \eqref{potential},
with $K$ replaced by $\eps/\tau$ with $\eps=1/|\bbD_0^{-1}|$ the 
positive-definiteness constant
of $\bbD_0$, i.e.\ $\bbD_0e\!\colon\!e\ge\eps|e|^2$. Then
\eqref{K-large} yields \eqref{tau-small}.
The strict convexity of the overall functional \eqref{potential} then follows
from the strict convexity of the $\varrho$- and $\kappa$-terms in
\eqref{potential}.
$\hfill\Box$

\medskip

Let us note that, as the physical dimension of $\bbC_1$ is Pa and of $\bbD_0$
is Pa\,s, the right-hand side of \eqref{tau-small} has the physical dimension
indeed seconds, as expected.

In particular, for sufficiently small time steps, the formula 
\eqref{IBVP-damage-small-disc}--\eqref{IBVP-damage-small-disc-BC}
possesses just one weak solution.

\begin{remark}[Insufficiency of viscosity in damage.]\label{rem-visco-dam}
  \upshape
  \NEW{It is notable that inventing some viscosity into the damage
    flow-rule (as sometimes considered in engineering literature
    \cite{KnRaBr01CPFA,MCMD07CDMC,MiWeHo10TCPF,WilKno09SDFO}) cannot
    stabilize the monolithic scheme. This is because } 
    the function 
$(e,\alpha)\mapsto\frac12\bbC(\alpha)e{:}e+\frac12K\alpha^2$ 
is not convex in general, no matter how $K$ is big. This can be seen even
for $d=1$ by analysing the positive definiteness of the corresponding
Hessian, whose determinant, in contrast to \eqref{det-hessian+},
would then contain the factor
$\frac12\gamma''(\alpha)\bbC_1^{}e{:} e
+K-\gamma'(\alpha)^2e\bbC_1^{}\bbC^{-1}(0)\bbC_1^{}e$
which, for $|e|\to\infty$, eliminates the influence of $K$ and
the positive definiteness is surely corrupted if $\bbC(0)$ is small (as usually
accepted as a modelling assumption).
\end{remark}

\begin{remark}[Fully implicit scheme.]\upshape
Note that the scheme 
\eqref{IBVP-damage-small-disc}--\eqref{IBVP-damage-small-disc-BC} is truly 
fully implicit only if $\taur=0$ because we took the viscous moduli
``delayed'', i.e.\ $\taur\bbC(\alpha_\tau^{k-1})$, instead of
$\taur\bbC(\alpha_\tau^{k})$ which would be a considerable alternative but
would corrupt the variational structure which was proved in
Proposition~\ref{prop} and which may be advantageous for some optimization
algorithms. Yet, e.g.\ the mentioned Newton-Raphson may not rely on
existence of some potentials and then such a fully implicit
scheme can be considered, even for $p>1$ since the last term in
\eqref{equil-in-w-disc} and its estimation \eqref{est-of-difference} below
would be avoided.
\end{remark}

\section{A-priori estimates and convergence for $\tau\to0$}\label{sec-anal}

Thorough the whole paper, 
we will use the standard notation for the Lebesgue $L^p$-spaces and
$W^{k,p}$ for Sobolev spaces whose $k$-th distributional derivatives 
are in $L^p$-spaces. We will also use the abbreviation $H^k=W^{k,2}$. 
In the vectorial case, we will write 
$L^p({\varOmega};\R^n)\cong L^p({\varOmega})^n$ 
and $W^{1,p}({\varOmega};\R^n)\cong W^{1,p}({\varOmega})^n$. 
On the time interval $I=[0,T]$, we consider the Bochner spaces 
$L^p(I;X)$ of Bochner measurable mappings $I\to X$ whose norm is in $L^p(I)$,
with $X$ being a Banach space, while $C_{\rm w}(I;X)$ will denote
the Banach space of weakly continuous mappings $I\to X$.
We recall using the notation $Q=I\times\Omega$ and $\Sigma=I\times\Gamma$.

Considering a fixed (sufficiently small to make \eqref{potential} convex) time 
step $\tau>0$ and $\{u_\tau^k\}_{k=0,...,T/\tau}$
with $T/\tau\in\mathbb N$, we define the piecewise-constant and the piecewise 
affine interpolants respectively by 
\begin{subequations}\label{def-of-interpolants}
\begin{align}\label{def-of-interpolants-}
&&&
\overline{u}_\tau(t)= u_\tau^k,\qquad\ \
\underline u_\tau(t)= u_\tau^{k-1},\qquad\ \
\underline{\overline u}_\tau(t)=\frac12u_\tau^k+\frac12u_\tau^{k-1},
&&\text{and}
&&
\\&&&\label{def-of-interpolants+}
u_\tau(t)=\frac{t-(k{-}1)\tau}\tau u_\tau^k
+\frac{k\tau-t}\tau u_\tau^{k-1}
&&\hspace*{-8em}\text{for }(k{-}1)\tau<t\le k\tau.
\end{align}\end{subequations}
For $\{(v_\tau^k,\alpha_\tau^k,f_\tau^k,g_\tau^k)\}_{k=0,...,T/\tau}$, the 
similar meaning has also $\underline{\overline v}_\tau$,
$\alpha_\tau$, $\overline\alpha_\tau$, $\overline f_\tau$, $\overline g_\tau$, etc. 
In terms of these interpolants, the scheme \eqref{IBVP-damage-small-disc}--\eqref{IBVP-damage-small-disc-BC} can be written ``more compactly'' as: 
\begin{subequations}\label{IBVP-damage-small-}
\begin{align}
&\DT u_\tau=\underline{\overline v}_\tau\ \ \text{ and }\ \ 
\varrho\DT v_\tau-{\rm div}\,\big(
\bbC(\overline\alpha_\tau)e(\overline u_\tau)
+\bbD(\underline\alpha_\tau)e(\underline{\overline v}_\tau)\big)
=\overline f_\tau\,,
\label{IBVP-damage--1}
\\
&\partial\zetai\big(\DT\alpha_\tau\big)
+\frac12\bbC'(\overline\alpha_\tau)e(\overline u_\tau){:}e(\overline u_\tau)
-{\rm div}(\kappa|\nabla\overline\alpha_\tau|^{p-2}\nabla\overline\alpha_\tau)
\ni\phi'(\overline\alpha_\tau)
\label{IBVP-damage--2}
\end{align}
with the boundary conditions 
\begin{align}\label{IBVP-damage-small-BC}
\Big(
\bbC(\overline\alpha_\tau)e(\overline u_\tau)
+\bbD(\underline\alpha_\tau)e\big(\underline{\overline v}_\tau\big)\Big)
\vec{n}=\overline g_\tau\quad\text{ and }\quad
\nabla\overline\alpha_\tau\!\cdot\!\vec{n}=0\quad\text{ on }\ \Sigma.
\end{align}\end{subequations}

\begin{proposition}[Discrete energetics.]\label{prop-4.1}
  Let the assumptions \eqref{model-ass} and \eqref{ass} 
  be fulfilled
  and $p>1$, and let
  $\tau_0>0$ be taken from \eqref{tau-small}.
Then the analog of \eqref{energetics} as an inequality
\begin{align}\nonumber
\mathscr{T}(\DT u_\tau(t))+\mathscr{E}(u_\tau(t),\alpha_\tau(t))
+\Big(1-\sqrt{\frac\tau{\tau_0}}\Big)
\!\int_0^t\!\!\varXi(\underline\alpha_\tau;\DT u_\tau,\DT\alpha_\tau)\,\d t
\qquad
\\[-.3em]\label{energetics-disc}\le
\mathscr{T}(v_0)+\mathscr{E}(u_0,\alpha_0)
+\!\int_0^t\!\!\mathscr{F}(\bar t_\tau,\DT u_\tau)\,\d t
\end{align}
holds for any $t=k\tau$ with $k=1,...,T/\tau$, with $\bar t_\tau(t):=k\tau$ for 
$t\in((k{-}1\tau,k\tau)$, provided  $\tau>0$ is sufficiently small
satisfying \eqref{tau-small}. 
\end{proposition}

\noindent{\it Proof.} To pursue the physical energy estimates, we
test the first equation in \eqref{IBVP-damage--1} by $\varrho\DT v_\tau$,
the second equation in \eqref{IBVP-damage--1} by $\DT u_\tau=\underline{\overline v}_\tau$, and \eqref{IBVP-damage--2} by $\DT\alpha_\tau$. The
first test leads, by the binomial formula, to the equality 
\begin{align}\nonumber\int_0^T\!\!\varrho\DT v_\tau{\cdot}\DT u_\tau\,\d t&=
\int_0^T\!\!\varrho\DT v_\tau{\cdot}\underline{\overline v}_\tau\,\d t=
\sum_{k=1}^{T/\tau}\varrho\frac{v_\tau^k{-}v_\tau^{k-1}}\tau{\cdot}
\frac{v_\tau^k{+}v_\tau^{k-1}}2
\\&=\sum_{k=1}^{T/\tau}\frac\varrho2|v_\tau^k|^2-\frac\varrho2|v_\tau^{k-1}|^2
=\frac\varrho2|v_\tau(T)|^2-\frac\varrho2|v_0|^2
\label{binomial}\end{align} 
a.e.\ on $\Omega$.
For the second and the third mentioned tests, 
we use the notation (\ref{energetics}b,c) to
write (\ref{IBVP-damage-small-disc}b,c)--\eqref{IBVP-damage-small-disc-BC}
shortly in an abstract way as
$$
\mathscr{T}'\frac{v_\tau^k{-}v_\tau^{k-1}}\tau
+\mathscr{E}'(u_\tau^k,\alpha_\tau^k)
+\mathscr{D}_{(v,\alpha)}'\Big(\alpha^{k-1};\frac{u_\tau^k{-}u_\tau^{k-1}}\tau,
\frac{\alpha_\tau^k{-}\alpha_\tau^{k-1}}\tau\Big)
\ni\mathscr{F}_\tau^k
$$
with $\mathscr{F}_\tau^k$ as in \eqref{energies-F} but with
$f_\tau^k$ and $g_\tau^k$ from
\eqref{IBVP-damage-disc-1}--\eqref{IBVP-damage-small-disc-BC2},
and we test it by $(u_\tau^k,\alpha_\tau^k){-}(u_\tau^{k-1},\alpha_\tau^{k-1})$.
The first term was already done in \eqref{binomial}. The 
second and the third terms can be estimated by 
using the semiconvexity proved in Lemma~\ref{lem}, cf.\ also
\cite[Remark 8.24]{Roub13NPDE} for this estimation technique.
Here, denoting $\mathscr{D}_0(v):=\int_\Omega\frac12\bbD_0e(v){:}e(v)\,\d x$
and $\varXi_0(v):=\int_\Omega\bbD_0e(v){:}e(v)\,\d x$
and taking some reference time step $\tau_0>0$,
it gives
\begin{align}\nonumber
  &\Big\langle
  \mathscr{D}_{(v,\alpha)}'\Big(\alpha^{k-1};\frac{u_\tau^k{-}u_\tau^{k-1}}\tau,
\frac{\alpha_\tau^k{-}\alpha_\tau^{k-1}}\tau\Big)+
\mathscr{E}'(u_\tau^k,\alpha_\tau^k),\frac{(u_\tau^k,\alpha_\tau^k){-}(u_\tau^{k-1},\alpha_\tau^{k-1})}\tau\Big\rangle
\\\nonumber
&\ge\varXi_0\Big(\frac{u_\tau^k{-}u_\tau^{k-1}}\tau\Big)
+\Big\langle\mathscr{E}'(u_\tau^k,\alpha_\tau^k),\frac{(u_\tau^k,\alpha_\tau^k){-}(u_\tau^{k-1},\alpha_\tau^{k-1})}\tau\Big\rangle
\\\nonumber
&=
\Big\langle\mathscr{D}_0'\Big(\frac{u_\tau^k{-}u_\tau^{k-1}}\tau\Big)+
\mathscr{E}'(u_\tau^k,\alpha_\tau^k),\frac{(u_\tau^k,\alpha_\tau^k){-}(u_\tau^{k-1},\alpha_\tau^{k-1})}\tau\Big\rangle
\\\nonumber
&=\Big\langle\frac1{\sqrt{\tau_0\tau}}\mathscr{D}_0'(u_\tau^k)
+\mathscr{E}'(u_\tau^k,\alpha_\tau^k),\frac{(u_\tau^k,\alpha_\tau^k){-}(u_\tau^{k-1},\alpha_\tau^{k-1})}\tau\Big\rangle
\\\nonumber
&\qquad-\Big\langle\frac1{\sqrt{\tau_0\tau}}\mathscr{D}_0'(u_\tau^{k-1})
,\frac{u_\tau^k{-}u_\tau^{k-1}}\tau\Big\rangle
  +\Big(1-\sqrt{\frac\tau{\tau_0}}\Big)\varXi_0\Big(\frac{u_\tau^k{-}u_\tau^{k-1}}\tau\Big)
\\&\nonumber\ge
\frac1\tau\bigg(\frac1{\sqrt{\tau_0\tau}}\mathscr{D}_0^{}(u_\tau^k)
+\mathscr{E}(u_\tau^k,\alpha_\tau^k)
-\frac1{\sqrt{\tau_0\tau}}\mathscr{D}_0^{}(u_\tau^{k-1})
-\mathscr{E}(u_\tau^{k-1},\alpha_\tau^{k-1})\bigg)
  \\\nonumber
&\qquad-\Big\langle\frac1{\sqrt{\tau_0\tau}}\mathscr{D}_0'(u_\tau^{k-1})
,\frac{u_\tau^k{-}u_\tau^{k-1}}\tau\Big\rangle
  +\Big(1{-}\sqrt{\frac\tau{\tau_0}}\Big)\varXi_0\Big(\frac{u_\tau^k{-}u_\tau^{k-1}}\tau\Big)
\\\label{est-by-semiconvex}
&=
\frac{\mathscr{E}(u_\tau^k,\alpha_\tau^k)-\mathscr{E}(u_\tau^{k-1},\alpha_\tau^{k-1})}
  \tau
  +\Big(1{-}\sqrt{\frac\tau{\tau_0}}\Big)\varXi_0\Big(\frac{u_\tau^k{-}u_\tau^{k-1}}\tau\Big)\,.
  \end{align} 
Here we use convexity of $\mathscr{E}+\mathscr{D}_0/\sqrt{\tau_0\tau}$
provided $\tau>0$ is small enough, which follows from the convexity
of $(e,\alpha)\mapsto\bbC(\alpha)e{:}e+\bbD_0e{:}e/\sqrt{\tau_0\tau}$
for $\tau\le\tau_0$ with $\tau_0$ from \eqref{tau-small}.

The resting contribution to the dissipation rate $\varXi$
is even with the factor 1. Altogether, 
\eqref{energetics-disc} is proved.
$\hfill\Box$

\medskip

\begin{proposition}[Numerical stability -- a-priori estimates.]\label{prop2}
Let  the assumptions of Proposition~\ref{prop-4.1} be fulfilled,
and let now $\tau\le\tau_0/2$ with $\tau_0$ from \eqref{tau-small}.
Then, the following a-priori estimates hold with $C$ independent of $\tau$:
\begin{subequations}\label{est}
\begin{align}\label{est1}
  &\|u_\tau\|_{H^1(I;H^1(\Omega;\R^d))}^{}\le C
\ \text{ \NEW{and} }\ 
\|u_\tau\|_{W^{1,\infty}(I;L^2(\Omega;\R^d))}^{}\le C,
\\\label{est2}
&\|v_\tau\|_{L^2(I;H^1(\Omega;\R^d))}^{}\le C
\ \ \text{ \NEW{and} }\ \  
\|u_\tau\|_{L^{\infty}(I;L^2(\Omega;\R^d))}^{}\le C,
\\&\label{est3}
\|\alpha_\tau\|_{L^\infty(I;W^{1,p}(\Omega))}^{}\le C
\ \ \text{ \NEW{and} }\ \ \ \|\alpha_\tau\|_{H^1(I;L^2(\Omega))}^{}\le C,
\intertext{and, if in addition $\varrho\in W^{1,p}(\Omega)$ with
  $p>2$ if $d=2$ or $p=3$ if $d=3$, then also}
&\label{est4}
\|\sqrt\varrho\DT v_\tau\|_{L^2(I;H^1(\Omega;\R^d)^*)}^{}\le C.
\end{align}\end{subequations}
\end{proposition}

\noindent{\it Proof.}
We have \eqref{energetics-disc} with the factor $(1{-}\sqrt{1/2})$
in front of the $\varXi$-term at disposal.
The right-hand side of \eqref{energetics-disc} estimated from above
by $\int_0^t C\|f(t)\|_{L^2(\Omega;\R^d)}^2+C\|g(t)\|_{L^2(\Omega;\R^d)}^2
+\frac14\varXi_0(\DT u_\tau(t))\,\d t$ with some $C$ large enough,
and then $\frac14\varXi_0(\DT u_\tau)$ can be absorbed in the dissipation
term in the left-hand side. From the coercivity of $\mathscr{T}$,
$\mathscr{E}$, and $\varXi(\alpha;\cdot,\cdot)$, 
we eventually obtain \eqref{est}.

The estimate \eqref{est4} can be proved by a comparison
$\varrho\DT v_\tau={\rm div}(\bbC(\overline\alpha_\tau)e(\overline u_\tau)
+\bbD(\underline\alpha_\tau)e(\underline{\overline v}_\tau))+\overline f_\tau$,
cf.\ \eqref{IBVP-damage--1}. More in detail, when testing it by
$z\in L^2(I;H^1(\Omega;\R^d))$ and using also the boundary conditions
\eqref{IBVP-damage-small-BC}, one arrives to
\begin{align}\nonumber
  \big\|\sqrt\varrho\DT v_\tau\big\|_{L^2(I;H^1(\Omega;\R^d)^*)}^{}
  &=\sup_{\|z\|_{L^2(I;H^1(\Omega;\R^d))}\le 1}\bigg(
  \int_\Sigma \overline g_\tau\cdot\frac{z}{\sqrt\varrho}\,\d S\d t
  \\&\nonumber
  \qquad\qquad+\int_Q\overline f_\tau\cdot\frac{z}{\sqrt\varrho}
-\big(\bbC(\overline\alpha_\tau)e(\overline u_\tau)
+\bbD(\underline\alpha_\tau)e(\underline{\overline v}_\tau)\big)\colon
\nabla\frac z{\sqrt\varrho}\,\d x\d t\bigg),
\end{align}
from which \eqref{est4} follows by the H\"older inequality and by using the
already obtained estimates (\ref{est}a,b), cf.\ also
\cite[Sect.\,6.4]{KruRou18MMCM}.
$\hfill\Box$

\medskip

The concept of weak solutions to the original continuous problem 
\eqref{KV-damage-} is a bit delicate because the damage driving force 
(or a ``pressure'') $\frac12\bbC'(\alpha)e(u)\!\colon\!e(u)$ is bounded only
in $L^1(Q)$ and testing it by $\DT\alpha\in L^2(Q)$
as occurs in \eqref{VI-ineq}
is not legitimate. To cope with this problem,
\NEW{this possibly not-integrable term} can be substituted by
\begin{align*}
\int_Q\frac12\bbC'(\alpha)e(u)\!\colon\!e(u)\DT\alpha\,\d x\d t&=
\int_\Omega\frac12\bbC(\alpha(T))e(u(T))\!\colon\!e(u(T))\,\d x
\\[-.4em]&\qquad-\int_Q\bbC(\alpha)e(u)\!\colon\!e(\DT u)\,\d x\d t
-\int_\Omega\frac12\bbC(\alpha)_0e(u_0)\!\colon\!e(u_0)\,\d x\,.
\end{align*}
Similar note is about the gradient term 
${\rm div}(\kappa|\nabla\alpha|^{p-2}\nabla\alpha)$.
This last term can be just integrated by parts, i.e.\ 
$\int_Q{\rm div}(\kappa|\nabla\alpha|^{p-2}\nabla\alpha)\DT\alpha\,\d x\d t$
$=$ $\int_\Omega\frac\kappa p|\nabla\alpha_0|^p-\frac\kappa p|\nabla\alpha(T)|^p
\,\d x$. Thus we obtain the following:

\begin{definition}[Weak solutions.]\label{def}
The \NEW{pair} $(u,\alpha)\in H^1(I;H^1(\Omega;\R^d))\times
C_{\rm w}(I;W^{1,p}(\Omega))$ is a weak solution to 
the initial-boundary-value problem \eqref{KV-damage} if
$\DT u\in C_{\rm w}(I;L^2(\Omega;\R^d))$, 
$
\zetai(\DT\alpha)\in L^1(Q)$, $0\le \alpha\le 1$ 
a.e.\ on $Q$, and 
\begin{subequations}\begin{align}\nonumber
&\forall v\in L^2(I;H^1(\Omega;\R^d))\,\cap\,H^1(I;L^2(\Omega;\R^d)),\ \ 
    v|_{t=T}^{}=0:
\\&\nonumber\hspace{1.5em}
\int_Q
\bbD(\alpha)e(\DT u)+\bbC(\alpha)e(u)\!\colon\! e(v)-\varrho\DT{u}\!\cdot\!\DT{v}\,\d x\d t
\\[-.5em]&\qquad\qquad\qquad\qquad\qquad\quad
=\int_\Omega\varrho v_0\!\cdot\! v(0,\cdot)\,\d x
+\int_Q\!f\!\cdot\! v\,\d x\d t+\int_{\Sigma}\!\!g\!\cdot\! v\,\d S\d t,
\label{7-very-weak-sln-damage}
\\\nonumber
&\!\!\forall z\in L^1(I;W^{1,p}(\Omega))\cap L^\infty(Q):\ \ 
\\&\hspace{1.5em}\nonumber\int_Q\!
\Big(
\frac12\bbC'(\alpha)e(u){:}e(u)
-\phi'(\alpha)\Big)z
+\kappa|\nabla\alpha|^{p-2}\nabla\alpha\!\cdot\!\nabla z+\zetai(z)\,\d x\d t
\\[-.3em]\nonumber
&\qquad\quad
\ge\!\int_Q
\zetai(\DT\alpha)+\bbC(\alpha)e(u)\!\colon\! e(\DT u)
\,\d x\d t
+\int_\Omega\bigg(
\frac12\bbC(\alpha(T))e(u(T)){:}e(u(T))-\phi(\alpha(T))
\\[-.3em]
&\qquad\qquad\qquad\quad
+\frac\kappa p|\nabla\alpha(T)|^p
-
\frac12\bbC(\alpha_0)e(u_0){:}e(u_0)
+\phi(\alpha_0)-\frac\kappa p|\nabla\alpha_0|^p
\bigg)\,\d x,
\label{id6-KV-damage}
\end{align}\end{subequations}
and if also the resting initial conditions 
hold, i.e.\ $u|_{t=0}=u_0$ and $\alpha|_{t=0}=\alpha_0$ {\rm(}while $\DT u|_{t=0}=v_0$
is already involved in \eqref{7-very-weak-sln-damage}{\rm)}.
\end{definition}

\begin{proposition}[Convergence for $\tau\to0$.]\label{prop3}
  Let the assumptions of Proposition~\ref{prop-4.1} be fulfilled now with
  $p>d$.
Then there exists a selected subsequence and its limit $(u,\alpha)$ such that
  \begin{subequations}\label{conv}
\begin{align}\label{conv1}
  &u_\tau\NEW{\rightharpoonup} u\ \text{ in }
  H^1(I;H^1(\Omega;\R^d))\ \NEW{(weakly)}
  \ \text{ and }\ 
  u_\tau\NEW{\stackrel{_*}{\rightharpoonup}} u\ \text{ in } L^\infty(I;L^2(\Omega;\R^d))\ \NEW{(weakly*)},
\\\label{conv2}
&v_\tau\NEW{\rightharpoonup} v
\ \,\text{ in }\ L^2(I;H^1(\Omega;\R^d))\ \NEW{(weakly)\ \ and\ \ }v=\DT u,
\\\label{conv3}
&\alpha_\tau\NEW{\stackrel{_*}{\rightharpoonup}}\alpha\ \text{ in }\ L^\infty(I;W^{1,p}(\Omega))\ \NEW{(weakly*)}
\ \text{ and }\ \alpha_\tau\NEW{\rightharpoonup}\alpha
\ \text{ in }\ H^1(I;L^2(\Omega))\ \NEW{(weakly)}
\intertext{for $\tau\to0$. Moreover, every such a subsequence converges also
  as}\label{conv4}
 &u_\tau\to u\ \text{in }\ H^1(I;H^1(\Omega;\R^d))\ \NEW{(strongly)},
  \end{align}\end{subequations}
  and every such a limit $(u,\alpha)$ is a weak solution to the
  initial-boundary-value problem \eqref{KV-damage} according Definition~\ref{def}.
\end{proposition}

\noindent{\it Proof.} 
After selecting a subsequence converging  weakly* in the topologies
indicated in \eqref{est} and using the 
Aubin-Lions theorem for the damage and then continuity of the
superposition operator induced by $\bbC(\cdot)$, we can 
pass to the limit first in the semilinear force-equlibrium
equation, obtaining (\ref{KV-damage}a,b,d,e) in the weak form
\eqref{7-very-weak-sln-damage}.

For the damage flow rule, we need the strong convergence \eqref{conv4},
however. The strategy for this strong convergence
is to use the new variable $w=u{+}\taur v$ as in \cite{LaOrSu10ESRM},
and denote the piecewise affine and the piecewise constant 
interpolants respectively as
\begin{align}
w_\tau=u_\tau+\taur v_\tau\ \ \ \text{ and }\ \ \ 
\overline w_\tau=\overline u_\tau+\taur\DT u_\tau=
\overline u_\tau+\taur\underline{\overline v}_\tau\,.
\end{align}
Then we can write the time-discrete approximation of the force equilibrium 
\eqref{IBVP-damage--1} equivalently as 
\begin{align}\label{equil-in-w-disc}
\frac{\varrho}{\taur}\DT w_\tau
-{\rm div}\big(\bbD_0e(\underline{\overline v}_\tau)
+\bbC(\overline\alpha_\tau)e(\overline w_\tau)\big)
=\overline f_\tau+\frac{\varrho}{\taur}\DT u_\tau+{\rm div}\big(
(\bbC(\overline\alpha_\tau)
{-}\bbC(\underline\alpha_\tau))e(\DT u_\tau)\big)
\end{align}
accompanied with the initial/boundary conditions (\ref{KV-damage}c,d)
and similarly also the mentioned
limit force equilibrium can be written in such a form, namely
\begin{align}\label{equil-in-w}
\frac{\varrho}{\taur}\DT w-{\rm div}\big(\bbD_0e(v)+\bbC(\alpha)e(w)\big)
=f+\frac{\varrho}{\taur}\DT u
\end{align}
accompanied with the initial/boundary conditions
$(\bbD_0e(v)+\bbC(\alpha)e(w))\vec{n}=g$ and $w(0)=w_0:=u_0+\taur v_0$. An 
important 
fact is that, by comparison, we also know $\DT w\in L^2(I;H^1(\Omega;\R^d)^*)$
because ${\rm div}(\bbD_0v)$ and ${\rm div}(\bbC(\alpha)e(w))$
belong to $L^2(I;H^1(\Omega;\R^d)^*)$ so that $\DT w$ and ${\rm div}(\bbD_0v)$ 
and ${\rm div}(\bbC(\alpha)e(w))$ are in duality with 
$w\in L^2(I;H^1(\Omega;\R^d))$. Thus, testing \eqref{equil-in-w} by $w$ and 
integrating over $Q$, we obtain the conservation of mechanical energy in the 
form
\begin{align}\nonumber
&\int_\Omega\frac{\varrho}{2\taur}|w(T)|^2+\frac12\bbD_0e(u(T)){:}e(u(T))\,\d x
+\int_Q\bbD_0\taur e(v){:}e(v)+\bbC(\alpha)e(w){:}e(w)\,\d x\d t
\\&=\int_\Omega\frac{\varrho}{2\taur}|w_0|^2+\frac12\bbD_0e(u_0){:}e(u_0)\,\d x
+\int_Q\Big(f{+}\frac{\varrho}{\taur}\DT u\Big){\cdot}w\,\d x\d t
+\int_\Sigma g{\cdot}w\,\d S\d t.
\label{cons-engr}
\end{align}

We further
test \eqref{equil-in-w} by $\overline w_\tau$, and integrate over the time 
interval $[0,k\tau]$ with $k=1,...,T/\tau$. We use 
\begin{align}\nonumber
&\int_Q\!\varrho\DT w_\tau\!\cdot\!\overline w_\tau\,\d x\d t=
\int_Q\!\varrho(\DT u_\tau+\taur\DT v_\tau)
\!\cdot\!(\overline u_\tau+\taur\underline{\overline v}_\tau)\,\d x\d t
\\&\nonumber=
\int_Q\!\varrho(\DT u_\tau+\taur\DT v_\tau)
\!\cdot\!(\underline{\overline u}_\tau+\taur\underline{\overline v}_\tau)
+(\DT u_\tau+\taur\DT v_\tau)\!\cdot\!
(\overline u_\tau-\underline{\overline u}_\tau)\,\d x\d t
\\&=\int_\Omega\frac\varrho2|u_\tau(T)+\taur v_\tau(T)|^2
-\frac\varrho2|u_0+\taur v_0|^2\,\d x
+\!\!\!\!\lineunder{\frac\tau2\int_Q\!\varrho(\DT u_\tau{+}\taur\DT v_\tau)\!\cdot\!
\DT u_\tau\,\d x\d t}{$=\mathscr{O}(\tau)$}
\label{dwdt-vs-w}\end{align}
where 
we used also \eqref{binomial}
and the estimate \eqref{est1} for the last integral. Further, we use  
\begin{align}\nonumber
&\int_0^T\bbD_0e(\underline{\overline v}_\tau){:}e(\overline w_\tau)\,\d t
=\int_0^T\!\!\taur\bbD_0e(\underline{\overline v}_\tau){:}e(\underline{\overline v}_\tau)
+\bbD_0e(\DT u_\tau){:}e(\overline u_\tau)\,\d t
\\[-.4em]&\ \ \ge\int_0^T\!\!\taur\bbD_0e(\underline{\overline v}_\tau)
{:}e(\underline{\overline v}_\tau)\,\d t+\frac12\bbD_0e(u_\tau(T)){:}e(u_\tau(T))
-\frac12\bbD_0e(u_0){:}e(u_0)
\end{align}
a.e.\ on $\Omega$.
We thus obtain the estimate
\begin{align}\nonumber
&\limsup_{\tau\to0}\int_Q\!\!\taur\bbD_0e(\underline{\overline v}_\tau)
{:}e(\underline{\overline v}_\tau)\,\d x\d t
\le\int_\Omega\frac\varrho{2\taur}|u_0{+}\taur v_0|^2+\frac12\bbD_0e(u_0){:}e(u_0)\,\d x
\\&\nonumber\qquad
-\liminf_{\tau\to0}\bigg(\int_Q
\bbC(\overline\alpha_\tau)e(\overline w_\tau){:}e(\overline w_\tau)\,\d x\d t
\\&\nonumber\qquad\qquad\quad
+\int_\Omega\frac\varrho{2\taur}|u_\tau(T){+}\taur v_\tau(T)|^2+\frac12\bbD_0e(u_\tau(T)){:}e(u_\tau(T))\,\d x\bigg)
\\&\nonumber\qquad+\lim_{\tau\to0}\bigg(\int_Q\overline f_\tau{\cdot}\overline w_\tau
+(\bbC(\underline\alpha_\tau){-}\bbC(\overline\alpha_\tau)e(\DT u_\tau)
\!\colon\! e(\overline w_\tau)
\,\d x\d t
+\int_\Sigma\overline g_\tau{\cdot}\overline w_\tau\,\d S\d t
+\mathscr{O}(\tau)\bigg)
\\&\nonumber\qquad
\le\int_\Omega\frac\varrho{2\taur}|u_0{+}\taur v_0|^2+\frac12\bbD_0e(u_0){:}e(u_0)\,\d x
+\int_Qf{\cdot}w-\bbC(\alpha)e(w){:}e(w)\,\d x\d t
\\&\nonumber\qquad\qquad\quad
-\int_\Omega\frac\varrho{2\taur}|u(T){+}\taur v(T)|^2
+\frac12\bbD_0e(u(T)){:}e(u(T))\,\d x+\int_\Sigma g{\cdot}w\,\d S\d t
\\&
\qquad
=\int_Q\taur\bbD_0e(v){:}e(v)\,\d x\d t,
\label{KV-visco-damageble-strong-Rothe}
\end{align}
where $\mathscr{O}(\tau)$ if from \eqref{dwdt-vs-w}. Here we used also
that $\sqrt\varrho v_\tau(T)\NEW{\rightharpoonup}\sqrt\varrho v(T)$
in $L^2(\Omega;\R^d)$, which follows from the second convergence in
\eqref{conv1} together with the estimate \eqref{est4}.
In \eqref{KV-visco-damageble-strong-Rothe}, we used also
\begin{align}\nonumber
&\bigg|\int_Q
\big(\bbC(\underline\alpha_\tau){-}\bbC(\overline\alpha_\tau)\big)
e(\DT u_\tau)\!\colon\! e(\overline w_\tau)\,\d x\d t\bigg|
\\&\qquad
\le\big\|\bbC(\underline\alpha_\tau){-}\bbC(\overline\alpha_\tau)
\big\|_{L^\infty(Q;\R^{d^4})}
\big\|e(\DT u_\tau)\big\|_{L^2(Q;\R^{d\times d})}
\big\|e(\overline w_\tau)\big\|_{L^2(Q;\R^{d\times d})}\to 0\,.
\label{est-of-difference}\end{align}
Here we used the first estimates in (\ref{est}a,b) together with 
the convergence \eqref{conv3} and the compact embedding of 
$L^\infty(I;W^{1,p}(\Omega))\,\cap\,H^1(I;L^2(\Omega))$ into $C(\overline Q)$ for 
$p>d$. This is actually 
the only spot where $p>d$ is vitally needed.
The last equality in \eqref{KV-visco-damageble-strong-Rothe} is the energy 
conservation in the mechanical-equilibrium part \eqref{cons-engr}.

As we already know $e(\underline{\overline v}_\tau)\NEW{\rightharpoonup}e(v)$
in $L^2(Q;\R^{d\times d})$, from \eqref{KV-visco-damageble-strong-Rothe}
we can see even the strong convergence \eqref{conv4}. 
Since $e(\DT u_\tau)=e(\underline{\overline v}_\tau)$,
it also says $e(\DT u_\tau)\to e(\DT u)$,
from which the desired strong
convergence $e(\overline u_\tau)\to e(u)$ needed for the limit passage in 
the damage flow rule follows. 

Now we can perform the limit passage in the discrete damage flow rule. 
This is, however, a bit technical. Actually, \eqref{IBVP-damage--2} means:
\begin{align}\nonumber
&\int_Q\zetai\big(z\big)
+
\frac12\bbC'(\overline\alpha_\tau)e(\overline u_\tau)
{:}e(\overline u_\tau)(z{-}\DT\alpha_\tau)
-\phi'(\overline\alpha_\tau)(z{-}\DT\alpha_\tau)
\\[-.5em]&\qquad\qquad\qquad\qquad
+\kappa|\nabla\overline\alpha_\tau|^{p-2}\nabla\overline\alpha_\tau
{\cdot}\nabla(z{-}\DT\alpha_\tau)\,\d x\d t
\ge
\int_Q\zetai\big(\DT\alpha_\tau\big)\,\d x\d t.
\label{IBVP-damage--2+}
\end{align}
Note that, on the time-discrete level, this inequality is indeed legitimate
since $\DT\alpha_\tau\in L^\infty(Q)$ so that it is in duality with 
$\frac12\bbC'(\overline\alpha_\tau)e(\overline u_\tau){:}e(\overline u_\tau)
\in L^1(Q)$; here we again use $p>d$ although some regularization could relax
this requirement at this point. Anyhow, $\DT\alpha_\tau$ is not 
uniformly bounded in $L^\infty(Q)$ and we must integrate this term 
by using also the discrete momentum equilibrium to lead 
\eqref{IBVP-damage--2+} closer to the limit inequality \eqref{id6-KV-damage}.
More specifically, we need to use a discrete analog of the identity 
$\frac12\DT\alpha\bbC'(\alpha)e(u)\!\colon\!e(u)=
\frac{\partial}{\partial t}\frac12\bbC(\alpha)e(u)\!\colon\!e(u)
-\bbC(\alpha)e(u)\!\colon\!e(\DT u)$ integrated over time
as used already in casting \eqref{id6-KV-damage} in Definition~\ref{def}.
Standardly, this gives the desired inequality if the 
$(e,\alpha)\mapsto\bbC(\alpha)e\!\colon\!e$ were convex. In our case,
we need to rely only on its (strengthened) semiconvexity:
actually, from \eqref{est-by-semiconvex}, we can also read the estimate
\begin{align}\nonumber
  &\frac{\mathscr{E}(u_\tau^k,\alpha_\tau^k)
    -\mathscr{E}(u_\tau^{k-1},\alpha_\tau^{k-1})}\tau\le
\Big\langle\mathscr{E}'(u_\tau^k,\alpha_\tau^k),\frac{(u_\tau^k,\alpha_\tau^k){-}(u_\tau^{k-1},\alpha_\tau^{k-1})}\tau\Big\rangle
+\sqrt{\frac\tau{\tau_0}}
\varXi_0\Big(\frac{u_\tau^k{-}u_\tau^{k-1}}\tau\Big)
\,.
\end{align} 
When summing it up for $k=1,...,T/\tau$ and writing in terms of the
interpolants, we have
\begin{align}\nonumber
&\int_\Omega\frac12\bbC(\alpha_\tau(T))e(u_\tau(T)){:}e(u_\tau(T))
-\frac12\bbC(\alpha_0)e(u_0){:}e(u_0)\,\d x
\\[-.4em]&\qquad\le
\int_Q\frac12\DT\alpha_\tau\bbC'(\overline\alpha_\tau)e(\overline u_\tau)
{:}e(\overline u_\tau)
+\bbC(\overline\alpha_\tau)e(\overline u_\tau){:}e(\DT u_\tau)
+\sqrt{\frac\tau{\tau_0}}\bbD_0e(\DT u_\tau){:}e(\DT u_\tau)\,\d x\d t\,.
\label{est-by-semiconvex++}\end{align}
Due to the estimate of $e(\DT u_\tau)$ in $L^2(Q;\R^{d\times d})$, the last term
\eqref{est-by-semiconvex++} converges to 0 as $\mathcal{O}(\sqrt\tau)$.

A similar note is about the term 
$\kappa|\nabla\overline\alpha_\tau|^{p-2}\nabla\overline\alpha_\tau
{\cdot}\nabla\DT\alpha_\tau$ whose treatment is however simpler, namely
$\int_\Omega\frac\kappa p|\nabla\alpha_\tau(T)|^p
-\frac\kappa p|\nabla\alpha_0|^p\d x\le\int_Q
\kappa|\nabla\overline\alpha_\tau|^{p-2}\nabla\overline\alpha_\tau
{\cdot}\nabla\DT\alpha_\tau\,\d x\d t$.
By all these estimates, \eqref{IBVP-damage--2+} yields also 
\begin{align}\nonumber
&\int_Q\zetai(z)
+\Big(
\frac12\bbC'(\overline\alpha_\tau)e(\overline u_\tau)
{:}e(\overline u_\tau)-\phi'(\overline\alpha_\tau)\Big)z
+\kappa|\nabla\overline\alpha_\tau|^{p-2}\nabla\overline\alpha_\tau
{\cdot}\nabla z\,\d x\d t
\\[-.5em]\nonumber
&\quad\ge
\int_Q\zetai\big(\DT\alpha_\tau\big)
+\bbC(\overline\alpha_\tau)e(\overline u_\tau){:}e(\DT u_\tau)\,\d x\d t
+\int_\Omega\bigg(\frac12\bbC(\alpha_\tau(T))e(u_\tau(T)){:}e(u_\tau(T))
\\[-.3em]
&\qquad\qquad
-\phi(\alpha_\tau(T))+\frac\kappa p|\nabla\alpha_\tau(T)|^p
-\frac12\bbC(\alpha_0)e(u_0){:}e(u_0)
+\phi(\alpha_0)-\frac\kappa p|\nabla\alpha_0|^p
\bigg)\,\d x\,.\end{align}
The limit passage towards the inequality \eqref{id6-KV-damage} is then
easy by weak (lower) semicontinuity, using also that
$e(u_\tau(T))\NEW{\rightharpoonup}e(u(T))$
in $L^p(\Omega;\R^{d\times d})$ and
$\nabla\alpha_\tau(T)\NEW{\rightharpoonup}\nabla\alpha(T)$
in $L^p(\Omega;\R^d)$.
$\hfill\Box$

\begin{remark}[\NEW{Convergence of iterative procedures.}]\upshape
\NEW{When applying a (here unspecified) conformal finite-element
space discretisation for the boundary-value problem
 \eqref{IBVP-damage-small-disc}--\eqref{IBVP-damage-small-disc-BC},
it becomes equivalent to a minimization of \eqref{potential} on the
respective finite-dimensional subspace $H^1(\Omega;\R^d)\times H^1(\Omega)$
respecting also the (discretised) constraints $0\le\alpha\le\alpha_\tau^{k-1}$.
This leads to a convex mathematical-programming problem.}
     The iterative procedures of the 
Newton-Raphson type (which are in the potential case then identical with 
sequential quadratic programming) enjoys a global guaranteed convergence
to the only one global minimizer \NEW{of \eqref{potential}. More in detail,
these methods need a globalization (so-called trust-region strategy) in general,
  but when the positive-definite Hessian has a bounded condition number,
  a global convergence is ensured provided the step-length used in
  particular iterations is controlled appropriately, cf.\ e.g.\ the monographs
  \cite{CoGoTo00TRM,NocWri06NO,SunYua06OTMN}. For this, we need here
  that the time-step criterion \eqref{tau-small} is satisfied. 
  Let us note that \eqref{potential} is then uniformly convex
  and the Hessian can effectively be considered 
  bounded. For this last attribute, one should consider the constraints on
  $\alpha$ valued in the interval $[0,1]$ and, because of the term
 $\gamma(\alpha)\bbC_1 e(u){:}e(u)$, also the constraints on $e(u)$ 
which is a-priori bounded in $L^2(\Omega;\R^{d\times d})$ uniformly due to 
the estimate \eqref{est3} and, when considering a fixed space discretisation,
also in $L^\infty(\Omega;\R^{d\times d})$. Thus the
condition number of the Hessian of \eqref{potential} is effectively bounded and 
iterative methods enjoy (at least theoretically) global convergence on the
mentioned trust region, although it may be very ill-conditioned
for fine space discretisations.
}
\end{remark}

\NEW{
  \begin{remark}[Uniqueness issue.]\upshape
    The uniquenees of the weak solution to the system
    \eqref{KV-damage}, which would in particular guarantee convergence
    of the whole sequence (not only selected subsequences) in
    Proposition~\ref{prop3}), seems open. Some results can be found e.g.\ in
    \cite[Prop.\,7.5.5]{KruRou18MMCM} but it requires $\bbD$ independent of
$\alpha$. On the other hand, a possible nonuniqueness is even a desired
phenomenon reflecting high sensitivity and unstability in real
experiments typically occuring in damage and fracture mechanics.
\end{remark}
  }

\section{Concluding remarks}\label{sec-rem}
Some special cases of the model and various modifications 
the model are worth mentioning:

\begin{remark}[Phase-field fracture.]\upshape
The concept of bulk damage can (approximately) imitate the philosophy of 
fracture along surfaces if the damage stored energy $\phi'$ is big.
The popular ansatz is 
\begin{align}
&\!\!\!\!\!\!
\varphi_{_{\rm E}}(e,\alpha,\nabla\alpha):=\gamma(\alpha)\bbC e{:}e
+\!\!\lineunder{G_{\rm c}\Big(\frac{1}{2\eps}(1{-}\alpha)^2\!
+\frac\eps p|\nabla\alpha|^p\Big)}{crack surface density}\!
\text{ with }\ \gamma(\alpha)=\frac{(\eps/\eps_0)^2{+}\alpha^2}2
\label{eq6:AM-engr+}
\end{align}
with $\eps_0>0$ and $G_{\rm c}>0$ fixed. For $p=2$, this
is known as the so-called \emph{Ambrosio-Tortorelli 
functional}. In the static case, this approximation was proposed in
\cite{AmbTor90AFDJ,AmbTor90AFDP} for the scalar case
and the asymptotic analysis for $\eps\to0$ was rigorously executed.
The generalization for the vectorial case is due to Focardi \cite{Foca01VAFD}.
Later, it was extended for evolution situation, namely for 
a rate-independent cohesive damage in \cite{Giac05ATAQ}, see also 
also \cite{BoFrMa08VAF,BoLaRi11TDMD,LaOrSu10ESRM,MieRou15RIST}
where also inertial forces are sometimes considered. It is vastly used
in engineering under the name a phase-field model, cf.\ e.g.\ 
\cite{MiWeHo10TCPF,SKBN18HAPF,TLBM18CNVP,WuNgu18LSIP} where also 
various modifications of the stored energy \eqref{eq6:AM-engr+} and various
dissipation potentials have been devised, although without any 
rigorous analytical justification of the resulted models.
In fact, this fits with
\eqref{energies-E} with 
$\phi(\alpha)=-G_{\rm c}(1{-}\alpha)^2/(2\eps)$ and $\kappa=\eps G_{\rm c}$.
The dissipation $\zeta=\zeta(\DT\alpha)$ is considered zero for
$\DT\alpha\le0$ and $+\infty$ otherwise. When accompanied by viscosity,
it fits with the assumptions \eqref{model-ass1} and \eqref{ass1}.
%
For a small length-scale parameter $\eps>0$, this model allows for a nearly 
complete damage and with a tendency to be localized on very small volumes 
along evolving surfaces where the cracks propagates, while elsewhere 
$\alpha\sim1$ not to make the term $(1{-}\alpha)^2/\eps$ 
too large. By this way, it approximates the infinitesimally thin cracks and
in particular their propagation. However, it is well known that such
phase-field model (as well as Griffith fracture model) does not work well for
initiation of cracks unless some notches are presented, cf.\ e.g.\ the
discussion in \cite{TLBM18CNVP}. Some attempt for improvement consists in 
nonquadratic degradation function $\gamma$, cf.\
\cite{BuOrSu13AFEA,SKBN18HAPF,WuNgu18LSIP}.
The monolithic discretisation is particularly suitable for such modification
because it does not rely on componentwise quadratic structure of
\eqref{eq6:AM-engr+} for $p=2$, in contrast to staggered schemes with
linear-quadratic solvers. On the other hand, to achieve the mentioned
good initiation and propagation and simultaneously allow for small
values of $\eps$, one should choose $\gamma$ such that
$\gamma'(1)=\mathcal{O}(1/\eps)$ but then $\tau_0$ in the criterion
\eqref{tau-small} is as $\mathcal{O}(\eps^2)$ for $\eps\to0$, which makes
this restriction quite strong.
\end{remark}
  
\begin{remark}[General $\phi$'s.]\upshape
  For a non-concave 
  $\phi$, one can modify
  \eqref{IBVP-damage-disc-2} by replacing $\phi'(\alpha_\tau^k)$ by the
  difference quotient $(\phi(\alpha_\tau^k){-}\phi(\alpha_\tau^{k-1}))/
  (\alpha_\tau^k{-}\alpha_\tau^{k-1})$ where-ever $\alpha_\tau^k\ne\alpha_\tau^{k-1}$.
  The approximate energetics \eqref{energetics-disc} as well
  as the a-priori estimates \eqref{est} keep holding, as well as the
  existence of a potential behind the incremental problems,
  cf.\ \cite{RouPan17ECTD}.
\end{remark}

\begin{remark}[Purely elastic models.]\label{rem-purely-elastic}\upshape
  The inertial term itself helps to convexify the incremental problems but,
  in contrast to the Kelvin-Voigt viscosity, this stabilization is
  dependent also on the spatial discretisation. Typically, the ratio
  between conditional numbers of the mass and the elasticity operator
  ${\rm div}(\bbC e(\cdot))$ after space discretisation is
  $\mathscr{O}(h^2)$ with $h>0$ the mesh size
  of the  spatial discretisation. This makes $\tau_0$ in the
  criterion \eqref{tau-small} with $\bbD_0$ replaced by the mass matrix
  dependent on $h$ as $\mathscr{O}(h^2)$ for $h\to0$, which would force 
  for practically unbearably too small time steps. On top of it, the convergence
  analysis of purely elastic model requires to enhance the model by
  strain gradients, cf.\ \cite[Sect.\,7.5.3]{KruRou18MMCM}, which is
  sometimes used in the context of damage mechanics to model dispersion of
  elastic waves \cite{BeRFAs12DRIG}. Moreover, sometimes even inertia gradients
  are used, see again \cite{BeRFAs12DRIG}, which then would have the same
  mesh-independent convexifying effect as the Kelvin-Voigt viscosity. 
\end{remark}
  

\begin{remark}[
    Plasticity, creep,
    \NEW{phase transformations}, 
    mass transfer.]
\upshape
Let us em\-p\-hasize that one can relatively routinely combine the above presented 
damage (or phase-field fracture) model with other phenomena when involving
some other internal variables as plastic/creep strain, \NEW{or
volume fraction of various phases (e.g.\ in shape-memory materials)},
or concentration \NEW{in poro-elastic materials}.
Then one immediately gets monolithic schemes for a combination
with plasticity (i.e.\ so-called ductile damage or fracture, in contrast
to brittle) or visco-elastic creep rheology, or for a combination
with 
a diffusant flow in damageable poro-visco-elastic media.
\end{remark}

\subsubsection*{Acknowledgment}
The author is deeply thankful for the hospitality of the University of
Seville during his stays in 2017-8 and for inspiring discussions with
Jos\'e Reinoso. 
\NEW{The discussion with Ladislav Luk\v san about optimization algorithms
  and the comments of an anonymous referee are highly appreciated.}
Also the partial support from the Czech Science Foundation through the grant
16-34894L ``Variational structures in continuum thermomechanics of solids'',
17-04301S ``Advanced mathematical methods for dissipative evolutionary 
systems'',
\NEW{18-03834S ``Localization phenomena in shape memory alloys'',}
and \NEW{19-04956S
``Dynamic and nonlinear behaviour of smart structures; modelling and
optimization''},
and also the institutional support RVO: 61388998 (\v CR) are
acknowledged.

\end{sloppypar}
\end{document}